\documentclass[12pt]{amsart}
\usepackage{amssymb, latexsym, amsthm}
\newtheorem{Thm}{Theorem}

\newtheorem{Cor}[Thm]{Corollary}

\newtheorem{Prop}[Thm]{Proposition}

\newcommand\adeles{{ad\`{e}les}}

\def\ie{{i.e.}}
\def\eg{{e.g.}}
\def\cf{{\it cf. \/}}

\def\lam{{\lambda}}

\def\s{\sigma}

\def\A{{\mathbb A}}

\def\F{{\mathbb F}}

\def\N{{\mathbb N}}

\def\R{{\mathbb {R}}}

\def\Z{{\mathbb Z}}

\def\B{{\mathcal B}}

\def\O{{\mathcal{O}}}

\def\Vals{{\mathcal{V}}}

\DeclareMathAlphabet{\mathscr}{OT1}{pzc}{m}{it}

\newcommand\tensor[1][{}]{{\otimes_{#1}}}

\def\ra{{\rightarrow}}

\def\minusset{{-}}

\def\sub{\subseteq}

\def\({\left(}
\def\){\right)}
\def\isom{{\;\cong\;}}

\def\normali{{\lhd}}

\newcommand\suchthat{{\,:\ \,}}
\newcommand\subjectto{{\,|\ }}

\newcommand\comp[1]{{{#1}^{\operatorname{c}}}}

\DeclareMathOperator{\Aut}{Aut}

\newcommand\Cayley[2]{{\operatorname{Cay}(#1;#2)}}

\DeclareMathOperator{\Gal}{Gal}

\newcommand\cond[2][!]{{\operatorname{cond}_{\if!#1\relax\else{\comp{#1}}\fi}(#2)}}

\newcommand\op[1]{{#1^{\operatorname{op}}}}

\renewcommand\L[2]{{\operatorname{L}^{#1}(#2)}}
\newcommand\M[1][d]{{\operatorname{M}_{#1}}}
\newcommand\GL[1][d]{{\operatorname{GL}_{#1}}}
\newcommand\PGL[1][d]{{\operatorname{PGL}_{#1}}}

\newcommand\PO[1][d]{{\operatorname{PO}_{#1}}}

\newcommand\PSL[1][d]{{\operatorname{PSL}_{#1}}}

\newcommand\abs[2][F]{|{#2}|_{#1}}

\newcommand{\set}[1]{{\{#1\}}}
\newcommand{\card}[1]{{\left|{#1}\right|}}
\newcommand\ideal[1]{{\left<{#1}\right>}}
\newcommand\sg[1]{{\ideal{#1}}}

\newcommand\dimcol[2]{{[{#1}\!:\!{#2}]}}

\newcommand\db[1]{{(\:\!\!({#1})\:\!\!)}}
\newcommand\mul[1]{{#1^{\times}}}
\newcommand\md[2][d]{{\mul{#2}/{\mul{#2}}^{#1}}}

\newcommand\valF{{\nu_0}}

\newcommand\Tref[1]{{Theorem \ref{#1}}}

\newcommand\Cref[1]{{Corollary \ref{#1}}}

\newcommand\eq[1]{{(\ref{#1})}}

\newcommand\Eq[1]{{Equation \eq{#1}}}

\newcommand\defin[1]{{\it{#1}}}
\long\def\half#1\halved{{\footnotesize{#1}}}
\long\def\forget#1\forgotten{}

\newcommand\dc[3]{{{#1}\backslash{#2}/{#3}}}
\newcommand\dom[2]{{{#1}\backslash{#2}}}

\newcommand\paper[6]{{{#1},\ {\it{#2}},\ {#3}\ {#4},\ {#5},\ ({#6}).}}
\newcommand\book[4]{{{#1},\ {{#2}},\ {#3},\ {#4}.}}

\newcommand\dd[4][!]{\abs[#2]{\:\!\det{\if!#1\relax\:\!\!\else(#1)\fi}}^{#3} #4}

\newcommand\binomq[3]{{\genfrac{[}{]}{0pt}{2}{#1}{#2}_{#3}}}
\newcommand\con[3]{{{#1}({#2},{#3})}}

\newif\ifXY
\XYfalse

\ifXY
\usepackage{xy}
\fi
\ifXY
\xyoption{all}
\fi

\begin{document}

\title
[Isospectral Cayley graphs]
{Isospectral Cayley graphs of some finite simple groups}

\def\HUJI{Inst. of Math., Hebrew Univ., Givat Ram, Jerusalem 91904,
Israel}
\def\BIU{Dept. of Math., Bar-Ilan University, Ramat-Gan 52900,
Israel}
\def\YALE{Dept. of Math., Yale University, 10 Hillhouse Av., New-Haven CT
06520, USA}

\author{Alexander Lubotzky}
\address{\HUJI}
\email{alexlub@math.huji.ac.il}

\author{Beth Samuels}
\address{\YALE
}
\email{beth.samuels@yale.edu}

\author{Uzi Vishne }
\address{\BIU}

\email{vishne@math.biu.ac.il}

\renewcommand{\subjclassname}{
      \textup{2000} Mathematics Subject Classification}

\date{Submitted: Dec. 30, 2004}

\maketitle

\begin{abstract}
We apply spectral analysis of quotients of the Bruhat-Tits
buildings of type $\tilde{A}_{d-1}$ to construct isospectral
non-isomorphic Cayley graphs of the finite simple groups
$\PSL[d](\F_q)$ for every $d \geq 5$ ($d \neq 6$) and prime power
$q > 2$.
\end{abstract}

\section{Introduction}

Let $X_1$ and $X_2$ be two finite graphs on $n$ vertices, and
$A_i$, $i = 1,2$, the $n\times n$ adjacency matrix of $X_i$. The
graphs $X_1$ and $X_2$ are \defin{isospectral} (or: cospectral),
if the multi-set of eigenvalues of $A_1$ is equal to that of
$A_2$. There are several methods to obtain isospectral graphs,
\eg\ Seidel switching method \cite{Ssm}, Sunada's method
\cite{Sunada}, \cite{Brooks}, \cite{Lu1} and various others.
However these methods produce isospectral graphs which are highly
non-homogeneous. Very little is known on isospectrality of Cayley
graphs. In fact, we are aware only of examples due to Babai
\cite{Babai}, who showed that the dihedral group of order $2p$
($p$ a prime) has at least $p/64$ pairs of generators which give
isospectral (non-isomorphic) Cayley graphs. Since dihedral groups
are almost abelian, the question remains for more complicated
groups. For example, can one find two sets of generators for the
symmetric group $S_n$, that will give isospectral non-isomorphic
Cayley graphs?

Here we prove the following
\begin{Thm}\label{main1}
For every $d \geq 5$ ($d \neq 6$) and every prime power $q$ and $e
\geq 1$ ($q^e > 4d^2+1$), there are two systems $A,B$ of
generators of the group $\PSL[d](\F_{q^{e}})$, such that the
Cayley graphs $\Cayley{\PSL[d](\F_{q^{e}})}{A}$ and
$\Cayley{\PSL[d](\F_{q^{e}})}{B}$ are isospectral and not
isomorphic.

The number of generators in $A$ and $B$ can be chosen to be either
$r = 2 \frac{q^d - 1}{q-1}$ or $r = \binomq{d}{1}{q} +
\binomq{d}{2}{q} + \cdots + \binomq{d}{d-1}{q}$ (where
$\binomq{d}{i}{q}$ denotes the number of subspaces of dimension
$i$ over $\F_q$ in the vector space $\F_q^d$).
\end{Thm}

Several remarks are in order here. First, our method produces
$\varphi(d)/2$ generating sets for which the Cayley graphs are
isomorphic (where $\varphi$ is Euler's function), which explains
the exceptional case $d = 6$. Secondly, if $A$ is a symmetric set
of generators for a group $G$, define the \defin{Cayley complex}
of $G$ with respect to $A$ to be the clique complex induced by the
Cayley graph (namely, $g_0,\dots,g_i \in G$ form an $i$-cell of
the complex iff each pair is connected in the Cayley graph). With
this definition, we prove that the Cayley complexes, rather than
graphs, are isospectral. Finally, for fixed $d$ and $q$, we obtain
infinitely many isospectral pairs which are $r$-regular with the
same $r$. This is of its own interest, since the usual methods to
construct isospectral $r$-regular graphs give pairs for which $r$
goes to infinity.

\begin{Thm}\label{main2}
Let $d \geq 3$ and let $r$ be as in \Tref{main1}. Then, for every
$m \in \N$, there exist $m$ isospectral non-isomorphic $r$-regular
graphs.
\end{Thm}

Using Sunada's method, Brooks \cite{Brooks} obtained such a result
for $r = 3$.

Let us now outline our method. The graphs considered in Theorems
\ref{main1} and \ref{main2} are in fact the $1$-skeletons of
Cayley complexes, or subgraphs of them.  The complexes are
obtained as quotients of the Bruhat-Tits building $\B_d(F)$
associated with the group $\PGL[d](F)$, where $F$ is a positive
characteristic local field. In spite of the `finite nature' of the
constructed objects, the proof is based on infinite dimensional
representation theory and the theory of division algebras over
global fields.

In more details, let $F = \F_q\db{t}$ be the field of Laurent
power series over a finite field $\F_q$ of order $q$, $G =
\PGL[d](F)$ and $\B_d(F)$ the Bruhat-Tits building associated with
$G$ (\cf \cite{Ronan},\cite{Tits}). The vertices of $\B = \B_d(F)$
can be identified with the cosets $G/K$ where $K = \PGL[d](\O)$ is
a maximal compact subgroup, and $\O = \F_q[[t]]$, the ring of
integers of $F$. Let $A$ be the adjacency operator or Laplacian
acting on $\L2{\B}$ (complex functions defined on the vertices).
$A$ commutes with the action of $G$, hence it is well defined on
$\dom{\Gamma}{\B}$ for every discrete subgroup $\Gamma$ of $G$. In
particular, if $\Gamma$ is a uniform lattice in $G$, \ie\ a
discrete cocompact subgroup, $A$ induces on $X = \dom{\Gamma}{\B}$
the adjacency matrix of the $1$-skeleton graph, which is an
$r$-regular graph for $r = \sum_{i=1}^{d-1}{\binomq{d}{i}{q}}$.

Now, $\L2{\dom{\Gamma}{G}}$ is a unitary representation space for
$G$. The following is well known (see \cite{Pesce} and also
\cite{paperI} or \cite{paperIII}):
\begin{Prop} \label{prop1}
Let $G = \PGL[d](F)$, $K =\PGL[d](\O)$ and let $\Gamma_1$ and
$\Gamma_2$ be two discrete cocompact subgroups of $G$. If
$\L2{\dom{\Gamma_1}{G}} \isom \L2{\dom{\Gamma_2}{G}}$ as (right)
$G$ representation spaces, then the simplicial complexes
$\dc{\Gamma_1}{G}{K}$ and $\dc{\Gamma_1}{G}{K}$ are isospectral.
In fact, they are even strongly isospectral, \ie\ with respect to
the higher dimensional Laplacians.
\end{Prop}

Using representation theory (Jacquet-Langlands correspondence) and
the theory of division algebras, we prepared in \cite{paperIII}
various examples of arithmetic subgroups of $\PGL[d](F)$ which
will provide the means to prove \Tref{main2}. This is a discrete
analog of Theorem 1 of \cite{paperIII}, claiming that for every $m
\in \N$ and $d \geq 3$, the symmetric space
$\PGL[d](\R)/\PO[d](\R)$ covers $m$ isospectral non-isomorphic
Riemannian manifolds. The proof when translated to our discrete
situation, gives $m$ isospectral simplicial complexes (details can
be found in \cite{paperIII}). Since in our situation the complex
is completely determined by its $1$-skeleton, we obtain $m$
$1$-skeletons which are not isomorphic as graphs.
It should be remarked that the proof of isospectrality in
\cite[Theorem~3]{paperIII} uses the Jacquet-Langlands
correspondence in positive characteristic; the details of the
proof of this correspondence appear in the literature for the
characteristic zero case (in general) and for prime $d$ in the
characteristic $p$ case. See Remark 1.6 in \cite{paperI} for more
discussion on this.

For the proof of \Tref{main1}, we use a family of remarkable
arithmetic lattices constructed by Cartwright and Steger
\cite{CS1}, which act simply transitively on the vertices of
$\B_d(F)$.  We denote their congruence subgroups by $\Delta(I)$,
where $I \normali \F_q[t]$ is a prime ideal. We proved in
\cite{paperII} that the quotient complexes
$\dom{\Delta(I)}{\B_d(F)}$ are Cayley complexes for a group which
lies between $\PGL[d](\F_{q^{e}})$ and $\PSL[d](\F_{q^{e}})$ (here
$e = \dim(\F_q[t]/I)$).  In particular, their $1$-skeletons are
Cayley graphs. We already have used similar complexes in
\cite{paperII} to give explicit constructions
of Ramanujan complexes  (see \cite{paperI}).
This time we summon them
to provide isospectral Cayley graphs.

The paper is organized as follows: In Section~\ref{sec:2} we
recall the notation and results from \cite{paperIII} and prove
\Tref{main2}. In Section~\ref{sec:3} we recall the
Cartwright-Steger lattices and their quotients, as described in
\cite{paperII}, and prove \Tref{main1}. In Section~\ref{explicit}
we present explicit examples of isospectral Cayley graphs.

\section*{Acknowledgements}

The authors acknowledge support of grants by the NSF and the
U.S.-Israel Binational Science Foundation.

\section{Proof of \Tref{main2}} \label{sec:2}
In \cite[Theorem~4]{paperIII}, we proved the following:
\begin{Thm}\label{fromIII}
Let $F$ be a local field of positive characteristic, $G =
\PGL[d](F)$ where $d \geq 3$, $K$ a maximal compact subgroup and
$\B_d(F) = G/K$ the associated Bruhat-Tits building. Then for
every $m \in \N$ there exists a family of $m$ torsion free
cocompact arithmetic lattices $\set{\Gamma_i}_{i = 1,\dots,m}$ in
$G$, such that the finite complexes $\dom{\Gamma_i}{\B_d(F)}$ are
isospectral and not commensurable.
\end{Thm}

The isospectrality of the $\dom{\Gamma_i}{\B_d(F)}$ proved in this
theorem is `representation theoretic isospectrality' i.e.
$\L2{\dom{\Gamma_i}{G}}$ are isomorphic as $G$-representation
spaces, where $G = \PGL[d](F)$. As explained there, this implies a
`strong combinatorial isospectrality' in the sense that not only
the spectrum of the standard Laplacians is independent of $i$, but
the same goes for the spectrum of the higher dimensional or
colored Laplacians.

We outline the construction here and refer the reader to
\cite{paperIII} for more details. Let $k$ be a global field
contained in $F$, and $\nu_0$ the valuation of $k$ with respect to
which the completion of $k$ is $F$. Since $F \isom \F_q\db{t}$,
one may choose $k = \F_q(t)$ and $\nu_0$ the $t$-adic valuation
(\ie\ $\nu_0(t^i) =i$ and $\nu(f) = 0$ for $f \in \F_q[t]$ prime
to $t$). The completion of $k$ with respect to a valuation $\nu$
is denoted $k_{\nu}$.

Now let $T$ be a finite set of valuations of $k$ not containing
$\nu_0$, and let $D_1,\dots,D_m$ be central division algebras of
degree $d$ over $k$, which are split by $k_{\nu}$ for every
valuation $\nu \not \in T$ (in particular by $F$), and remain
division algebras over $k_{\theta}$ for every $\theta \in T$.
Moreover we require that $D_i$ is not isomorphic as a ring to
either $D_j$ or $\op{D_j}$ for $i \neq j$ (recall that division
algebras over $k$ are uniquely determined---as algebras---by their
local invariants, which are in $\set{0,1/d,\dots,(d-1)/d}$).   We
showed in \cite{paperIII} that $T$ can be chosen to accommodate
$m$ such different division rings if $2^{\card{T}}/\card{T}
> 2m$.

Let $G'_i$ denote the algebraic group defined as the
multiplicative group  $\mul{D_i}$ modulo its center. Our
assumptions guarantee that $G'_i(k_{\theta})$ is a compact group
for every $\theta \in T$, and that $G'_i(F) = \PGL[d](F)$. Define
the ring
\begin{equation}\label{R0Tdef}
R_{0,T} = \set{x \in k \suchthat \nu(x)\geq 0 \quad \mbox{for
every $\nu \not\in T \cup \set{\valF}$}}.\end{equation} In
\cite[Section~2]{paperIII} we give an explicit description of
$G'_i$ as an algebraic group defined over $R_{0,T}$. Since
$R_{0,T}$ is discrete in the product $\prod_{\nu \in T \cup
\set{\nu_0}} {k_{\nu}}$, it follows that $G'_i(R_{0,T})$ are well
defined discrete subgroups of $\PGL[d](F)$.

Fix an ideal $I \normali R_{0,T}$, and let $\Gamma_i =
\con{G'_i}{R_{0,T}}{I}$ be the congruence subgroup (namely the
kernel of the map $G'_i(R_{0,T}) \ra G'_i(R_{0,T}/I)$ induced by
the projection $R_{0,T} \ra R_{0,T}/I$).

We prove \Tref{fromIII} above in \cite{paperIII}, by showing that
$\dom{\Gamma_i}{\B_d(F)}$ are isospectral and non-commensurable
(namely they have no finite sheeted common cover). To explain how
this result implies \Tref{main2}, we briefly describe the
Bruhat-Tits buildings associated with $G = \PGL[d](F)$. Let $\O
\sub F$ be the valuation ring, then $K =
\PGL[d](\O)$ is a maximal compact subgroup of $G$.

Consider the set of $\O$-submodules of $F^d$ of maximal rank (such
as $\O^d$), modulo the equivalence relation $M \sim cM$ for every
$c \in \mul{F}$. The action of $\GL[d](F)$ on bases of $F^d$
induces a transitive action of $G$ on the equivalence classes of
modules, with $K$ being the stabilizer of $[\O^d]$. Therefore, we
can identify classes of submodules with cosets in $G/K$. To define
the graph structure, we take the classes of modules as vertices,
and connect two classes $\chi,\chi'$ by an edge iff there are
representatives $M \in \chi$ and $M' \in \chi'$ such that $t M
\subset  M' \subset M$ (the set of neighbors of a vertex $[M] \in
\B^{(0)}$ is in one to one correspondence with the set of vector
subspaces of $M/tM \isom \F_q^d$ over $\F_q$, and so the valency
is $\binomq{d}{1}{q} + \binomq{d}{2}{q} + \dots +
\binomq{d}{d-1}{q}$). The resulting graph is the $1$-skeleton of
the building.

For the higher dimensional structure we take the clique complex
defined by the $1$-skeleton, namely the $i$-cells are the complete
subgraphs of size $i+1$. Equivalently, $\chi_0,\dots,\chi_i$ form
an $i$-cell iff there are representatives $M_j \in \chi_j$ such
that $t M_0 \subset M_i \subset M_{i-1} \subset \dots \subset M_1
\subset M_0$, after rearrangemenet. In particular there are no
$d$-cells. We let $\B_d(F)$ denote the resulting complex.

If $\Gamma \sub G$ is a cocompact lattice, then the quotient
$\dom{\Gamma}{\B_d(F)}$ is a finite complex (the complex structure
inherited from $\B_d(F)$). In order for the projection $\B_d(F)
\ra \dom{\Gamma}{\B_d(F)}$ to be a local isomorphism, $\Gamma$ has
to be torsion free. However a cocompact lattice of $G$ always has
a finite index torsion-free subgroup, and for arithmetic groups we
may choose the ideal $I$ small enough so that the $\Gamma_i$ are
all torsion free.

The complexes $\dom{\Gamma_i}{\B_d(F)}$ are isospectral by
\Tref{fromIII}, so the underlying graphs are also isospectral. On
the other hand if two of these graphs are isomorphic, then they
define isomorphic clique complexes, which again contradicts
\Tref{fromIII}. This completes the proof of \Tref{main2} for $r =
\binomq{d}{1}{q} + \binomq{d}{2}{q} + \cdots +
\binomq{d}{d-1}{q}$.

The determinant $G/K \ra \Z/d$ induces a $d$-coloring of the
(directed) edges: if $t M \sub M' \sub M$, the edge from $M$ to
$M'$ has color $\dim_{\F_q}(M/M')$. It is easily seen that all
edges connected to $[M]$ within one cell, have distinct colors.
The colored Laplacian $A_\ell$, $\ell = 1,\dots,d-1$, is then
defined on $\L2{G/K}$ by summing over neighbors of color $\ell$,
namely $(A_{\ell}f)(\chi)$ is the sum of $f(\chi')$ over all
neighbors of $\chi$ for which $(\chi,\chi')$ has color $\ell$. Of
course, $A_1+\dots+A_{d-1}$ is the standard Laplacian $A$ of
$\B_d(F)$. If two quotients $\dom{\Gamma_1}{\B_d(F)}$ and
$\dom{\Gamma_1}{\B_d(F)}$ of $\B_d(F)$ are
representation-theoretic isospectral, namely
$\L2{\dom{\Gamma_1}{G}} \isom \L2{\dom{\Gamma_2}{G}}$ as
$G$-spaces, then they are also isospectral in terms of $A_\ell$
for every $\ell$ (see \cite[Section~2]{paperI} for details).

Now let $\B^1$ denote the subgraph of $\B_d(F)$ defined on all
vertices, with edges only of color $1$. From the quotients
$\dom{\Gamma_i}{\B^1}$ we obtain $r$-regular graphs for $r =
2\frac{q^d-1}{q-1}$, which are isospectral with respect to the
standard Laplacian, in this case $A_1 + A_{d-1}$. On the other
hand, in \cite[Proposition~2.3]{paperII} we show that the color
$1$ part of the graph determines the whole skeleton. Therefore if
some quotients of
$\B^{1}$ are isomorphic, then so are the respective quotients of
$\B_d(F)$ which is impossible, by \Tref{fromIII}.  This completes
the proof of \Tref{main2}.

\section{Proof of \Tref{main1}} \label{sec:3}

In order to prove \Tref{main1} we describe remarkable lattices
$\Delta$ in $\PGL[d](F)$, which act simply transitively on the
vertices of $\B_d(F)$. The building can then be identified with
the Cayley complex of $\Delta$ with respect to a certain set of
generators (where the Cayley complex is the clique complex defined
by the Cayley graph). These lattice were introduced in \cite{CS1},
and used in \cite{paperII} to give explicit Ramanujan complexes.

We are given a prime power $q$ and an integer $d \geq 2$. The
local field $F = \F_q\db{t}$ contains $k = \F_q(t)$. Let $\nu_0$
denote the $t$-adic valuation, in which $F$ is complete, and let
$\O = \F_q[[t]]$ be the valuation ring. Let $\s$ be a fixed
generator of the Galois group of $\F_{q^d}/\F_q$, extended to
$\F_{q^d}(t)$ by acting trivially on $t$, and let
\begin{equation}\label{Ddef}
D^{(\s)} = \F_{q^d}(t)[z \subjectto z a = \s(a)z, \, z^d = 1+t],
\end{equation}
be a central division algebra of degree $d$ over $k$.

Let $F_1 = F \tensor[\F_q] \F_{q^d} = \F_{q^d}\db{t}$, and note
that $F_1/F$ is unramified in any valuation of degree $1$ of $k$
(\ie\ the minus degree valuation, denoted henceforth by
$\nu_{1/t}$, or a valuation induced by a linear prime $p \in
\F_q[t]$, denoted by $\nu_p$). Since $\nu_0(1+t) = 0$, local class
field theory guarantees that $1+t$ is a norm in $F_1/F$
\cite[Chapter~17]{Pierce}. But $D^{(\s)} \tensor[k] F = F_1[z]$
with $z^d = 1+t$, so by Wedderburn's criterion
\cite[Corollary~1.7.5]{Jac}, $D^{(\s)}$ splits over $F$. There is
therefore a natural embedding $D^{(\s)} \hookrightarrow \M[d](F)$,
which we will describe in detail later. As in the previous
section, $D^{(\s)}$ gives rise to the algebraic group
${G'}^{(\s)}$, defined as the multiplicative group of $D^{(\s)}$
modulo its center. Then we have ${G'}^{(\s)}(k) \hookrightarrow
{G'}^{(\s)}(F) = \PGL[d](F)$.

Let $\Omega$ denote the set of elements $u (1 - z^{-1}) u^{-1} \in
{G'}^{(\s)}(k)$, where $u$ runs over $\mul{\F_{q^d}}/\mul{\F_q}$,
and let $\Delta^{(\s)}$ be the subgroup of ${G'}^{(\s)}(k)$
generated by $\Omega$. It is shown in
\cite[Proposition~4.9]{paperII} that the embedding ${G'}^{(\s)}(k)
\hookrightarrow \PGL[d](F)$ takes $\Delta^{(\s)}$ into
$\PGL[d](R_0)$ for
\begin{equation}\label{R0def}
R_0 =
\F_q[1/t]
\end{equation}
(and the index $\dimcol{{G'}^{(\s)}(R_0)}{\Delta^{(\s)}}$ only
depends on $q$ and $d$, Proposition 3.5 there). Recall that $K =
\PGL[d](\O)$ is a maximal compact subgroup of $\PGL[d](F)$. In
\cite[Proposition~4.8]{paperII} we show (following \cite{CS1})
that $\Delta^{(\s)}$ acts simply transitively on the vertices
$\PGL[d](F)/K$ of the Bruhat-Tits building, or equivalently that
$$\Delta^{(\s)} \cdot K = \PGL[d](F) \qquad \mbox{and} \qquad \Delta^{(\s)} \cap K = 1.$$
Therefore we can identify (the vertices of) $\B^{0} = G/K$ with
(the elements of) $\Delta^{(\s)}$: every coset $g K$ has a unique
representative $\delta K$ with $\delta \in \Delta^{(\s)}$.

The complex structure can also be recovered from $\Delta^{(\s)}$,
as follows. The reduced norm of $1-z^{-1}$ and its conjugates is
$t/(1+t)$ \cite[Proposition~4.1]{paperII}, which is equivalent
modulo $\mul{\O}$ to the uniformizer $t$. It then follows that the
neighbors of color $1$ of $[\O^d]$ are $[\omega \O^d]$ for $\omega
\in \Omega$. Define $\hat{\Omega}$ to be the set of products
$\omega_1\dots \omega_\ell$, where $\omega_1,\dots,\omega_\ell \in
\Omega$, for which there exists $\omega_{\ell+1},\dots,\omega_d
\in \Omega$ such that the product
$\omega_1\dots\omega_\ell\omega_{\ell+1}\dots \omega_d$ equals $1$
(in ${G'}^{(\s)}(k)$). In Proposition 2.3 and Sections 5 and 6 of
\cite{paperII} we show that the neighbors (of arbitrary color) of
$[g \O^d]$ are $[g \omega \O^d]$ for $\omega \in \hat{\Omega}$,
and that $\B_d(F)$ is the Cayley complex of the group
$\Delta^{(\s)}$ with respect to the generators $\hat{\Omega}$.

The algebras $D^{(\s)}$ ramify at exactly two places, $T =
\set{\nu_{1/t},\nu_{1+t}}$. The ring defined by \Eq{R0Tdef} in
this situation is $R = R_{0,T} = \F_q[t, 1/t, 1/(1+t)]$. Let $I$
be an ideal of $R$, and set $L = R/I$. Since $\Delta^{(\s)} \sub
{G'}^{(\s)}(R_0) \sub {G'}^{(\s)}(R)$, we can denote by
$\Delta^{(\s)}(I)$ the kernel of the map $\Delta^{(\s)} \ra
\PGL[d](R/I)$ induced by the projection $\PGL[d](R) \ra
\PGL[d](R/I)$. In \cite[Theorem~6.6]{paperII} we prove that
$\Delta^{(\s)}_I = \Delta^{(\s)}/\Delta^{(\s)}(I)$ is a subgroup
of $\PGL[d](L)$ which contains $\PSL[d](L)$. Moreover, the index
of $\PSL[d](L)$ in $\Delta^{(\s)}_I$ is equal to the order of
$t/(1+t)$ in the group $\md[d]{L}$
(\cite[Proposition~6.7]{paperII}), and therefore is independent of
$\s$. Choosing the ideal $I$ properly, we can guarantee that
$\Delta^{(\s)}_I = \PSL[d](L)$ (Section~7 of \cite{paperII}). The
explicit embedding of $\Delta^{(\s)}_I$ into $\PGL[d](L)$ is given
in \cite[Section~9]{paperII}. We elaborate on this in the next
section.

To conclude, the set $\hat{\Omega}$ (viewed as elements of
$\PGL[d](L)$) generate a group $\Delta^{(\s)}_I$ which resides
between $\PSL[d](L)$ and $\PGL[d](L)$, and the resulting Cayley
complex is the quotient $\dom{\Delta^{(\s)}(I)}{\B_d(F)}$.

So far we merely described the Cayley complexes involved. To
finish the proof of \Tref{main1}, we quote again results proved in
\cite{paperIII}. Let $D_1$ and $D_2$ be two central division
algebras of degree $d$ over $k$, with the same set $T$ of
ramification places, but with different invariants at each
ramified place (and $\nu_0 \not \in T$). Define the algebraic
groups $G'_1$ and $G'_2$ as in Section~\ref{sec:2}. For simplicity
(and since this is all we need here), assume $\card{T} = 2$. In
\cite[Theorem~9]{paperIII} we proved that if $D_2$ is isomorphic
as a ring neither to $D_1$ nor to $\op{D_1}$, then the complexes
$\dom{\Gamma_1}{\B_d(F)}$ and $\dom{\Gamma_2}{\B_d(F)}$ cannot be
isomorphic for any finite index torsion-free subgroups $\Gamma_i
\sub G'_i(R_{0,T})$.

The invariants of $D^{(\s)}$ defined in \Eq{Ddef} are
$\frac{-t}{d},\frac{t}{d}$ where $t$ is chosen so that $\s^t$ is
the Frobenius automorphism on $\F_{q^d}$. Assume $\s_1,\s_2$ are
two generators of $\Aut(\F_{q^d}/\F_q)$ such that $\s_2 \neq
\s_1,\s_1^{-1}$. This situation is possible once $\phi(d) > 2$,
namely when $d \geq 5$, $d \neq 6$. Since the local invariants of
the algebras $D_i = D^{(\s_i)}$ defined in \Eq{Ddef} are
different, Proposition 11 of \cite{paperIII} guarantees that $D_2$
is not isomorphic (as a
ring) to $D_1$, nor to $\op{D_1}$.
Since we already know that $\B_d(F) \isom \Delta^{(\s_i)}$, the
above result implies that $\Delta^{(\s_1)}_I$ and
$\Delta^{(\s_2)}_I$ cannot be isomorphic to each other, when
viewed as Cayley complexes with the respective sets of generators
$\hat{\Omega}$. But the Cayley complex is completely determined by
the Cayley graph, and so $\Delta^{(s_1)}_I$ and $\Delta^{(s_2)}_I$
cannot be isomorphic as graphs.

In fact, the results of \cite{paperIII} show that \emph{as
complexes}, $\Delta^{(s_1)}_I$ and $\Delta^{(s_2)}_I$ are also
non-commensurable. Of course, they are commensurable as graphs by
Leighton's result \cite{Leighton} that every two graphs with the
same universal cover have a finite sheeted common cover.

It remains to show that the $\Delta^{(\s_i)}_I$ are isospectral to
each other. Let $\Vals$ denote the set of valuations of $\F_q(t)$,
and let $\A$ be ring of \adeles\ over $k$, namely the restricted
product of the completions $k_{\nu}$. Also let $\A_{\comp{T}}$ be
the produce over $\Vals \minusset T$, and notice that
$G'_i(\A_{\comp{T}}) = G({\comp{T}})$ for $i = 1,2$, since both
$D_i$ are split by $k_{\nu}$ for $\nu \in \Vals \minusset T$.

For a valuation $\nu \in \Vals$, Let $\O_{\nu}$ denote the
valuation ring in the completion $k_{\nu}$, and let $P_{\nu}$
denote the valuation ideal. Set
\begin{equation}\label{rdef}
r_{\nu} = \min\set{\nu(a)\suchthat a\in I},
\end{equation}
and let
\begin{equation}\label{Urdef}
U^{(r)} = \prod_{\nu \in \Vals - (T\cup \set{\nu_0})} \con{G}{\O_{\nu}}{P_{\nu}^{r_{\nu}}}.
\end{equation}
an open compact subgroup of $G(\A_{\comp{T}\minusset \nu_0})$.

Let $J_i^0$ be the product of $G'_i(k_{\nu_{1+t}})$ and the
pre-image of the standard unipotent subgroup of $\GL[d](\F_{q})$
under the map $G'_i(k_{\nu_{1/t}}) \ra \GL[d](\F_q)$ defined by
taking the matrices modulo $1/t$. By Definition 4.7 and
Proposition 4.8 in \cite{paperII}, $\Delta^{\s_i} = G'_i(k) \cap
G(\A_{\comp{T}}) J_i^{0}$, and so
$$\Delta^{\s_i} = G'_i(k) \cap G(F) U^{(r)} J_i^{0}.$$

Since $\dimcol{G'_i(k_{\nu_{1+t}})G'_i(k_{\nu_{1/t}})}{J_i^{0}}$
is independent of $i$, Proposition~15 of \cite{paperIII} implies
that
\begin{equation*}
\L2{\dom{\Delta^{\s_1}(I)}{G(F)}} \isom
\L2{\dom{\Delta^{\s_2}(I)}{G(F)}}
\end{equation*}
as $G(F)$-spaces. This implies isospectrality of the complexes
$$\dc{\Delta^{\s_i}(I)}{G(F)}{G(\O)} \isom
\dom{\Delta^{\s_i}(I)}{\B_d(F)}$$ for $i = 1,2$, by
\cite[Proposition~2]{paperIII}. But $\B_d(F) \isom
\Delta^{(\s_i)}$, so we have the isospectrality of
 $\Delta^{(\s_i)} /
\Delta^{(\s_i)}(I)$ (again, as Cayley complexes with respect to
the generators $\hat{\Omega}$). Choosing $I$ suitably (see
\cite[Section~7]{paperII}), we can assume $\Delta^{(\s_i)} /
\Delta^{(\s_i)}(I) \isom \PSL[d](R/I)$, as asserted.

\section{Explicit Examples of Isospectral Cayley Graphs}
\label{explicit}

Let $d \geq 5$ ($d \neq 6$) and let $q$ be a given prime power. In
Algorithm 9.2 of \cite{paperII} we give an explicit construction
of Ramanujan Cayley complexes of $\PSL[d](\F_{q^e})$ (see
\cite[Theorem~7.1]{paperII}), which is bound to work if $q^e
> 4d^2+1$. (This assumption is only used in
\cite[Proposition~7.3]{paperII}, and can often be ignored. At any
rate if $q = 2$ then we must take $e > 1$ since
$\F_2[t,1/t,1/(1+t)]$ has no quotients of dimension $1$).
Repeating the construction with generators $\s_1,\s_2$ of
$\Gal(\F_{q^d}/\F_q)$ such that $\s_2 \neq \s_1,\s_1^{-1}$, this
realizes the proof of \Tref{main1}, and provides $\phi(d)/2$
isospectral non-isomorphic Cayley complexes for any group
$\PSL[d](\F_{q^e})$. We illustrate this construction here with $d
= 5$, $q = 3$ and $e = 1$ (which is the smallest case possible).

\newcommand\MatrixFive[5]{{\(\!\!\! \begin{array}{ccccc}#1 \\[-0.0cm] #2 \\[-0.0cm] #3 \\[-0.0cm] #4 \\[-0.0cm] #5\end{array}\!\!\! \)\!
}}
\newcommand\LineFive[5]{{#1}\!\! &\!\! {#2}\!\! &\!\! {#3}\!\! &\!\! {#4}\!\! &\!\!
{#5}}

In order to view $\F_{q^d} = \F_{243}$ explicitly, take the
irreducible polynomial $\lam^5-\lam-1$ over $\F_q = \F_3$. We may
write $\F_{243} = \F_3[t \subjectto t^5 = t+1]$. We fix the
ordered basis $\set{1,t,t^2,t^3,t^4}$. {}From now on a linear
transformation of $\F_{243}$ (as a vector space over $\F_3$) is
represented via the chosen basis by the $5\times 5$ matrix. In
particular $\mul{\F_{243}}$ embeds in $\GL[5](\F_3)$ by the
regular representation, and the Frobenius automorphism $\phi$,
namely exponentiation by $3$, is represented by
$$\varphi_1 = \MatrixFive
{\LineFive{1}{0}{0}{1}{0}}
{\LineFive{0}{0}{1}{1}{0}}
{\LineFive{0}{0}{1}{0}{1}}
{\LineFive{0}{1}{0}{0}{2}}
{\LineFive{0}{0}{0}{1}{1}}.$$ We set $\s_1 = \phi$ and $\s_2 =
\phi^2$. The representing matrices are $\varphi_1$ and $\varphi_2
= \varphi_1^2$, respectively.

Define the polynomial $y(\lam) = (1+\lam)^5 - 1$ (this amounts to
taking $\beta = 1$ in Step (1) of \cite[Algorithm~9.2]{paperII}).
Take $\alpha = 1$, so that $\gamma = y(1) = 2^5 - 1 = 1$ in Step
(2) of the Algorithm. The minimal polynomials of these elements
are $p(\lam) = g(\lam) = \lam - 1$. (Here $g(t)$ must not be
invertible in $R_{0,T}$, namely $g(\lam) \neq \lam, \lam+1$). It
is possible to construct the Cayley complexes over the local ring
$\F_3[t]/\ideal{g(t)^s}$ for every $s \geq 1$. We continue with $s
= 1$, so in the notation of Step (3) in the Algorithm, $\bar{L} =
\F_3$ and $x = 1$. For $i = 1,2$, let $z_i$ be the matrix $(1 +
\beta x) \varphi_i = (1 + 1 \cdot 1) \varphi_i = 2 \varphi_i = -
\varphi_i$ (since $\varphi_i a = a^{3^i} \varphi_i$ for $a \in
\F_{243}$, conjugation by $z_i$ induces the automorphism $\s_i$,
as in the definition of $D_i$ in \Eq{Ddef}. Moreover $z_i^5 = -
\varphi_i^5 = -1 = 1+\alpha$).

Now we define the elements of $\Omega^{(\s_i)}$. Let $b^{(i)} = 1
- z_i^{-1} = 1 + \varphi_i^4$, namely
$$b^{(1)} = \MatrixFive
{\LineFive{2}{1}{2}{0}{1}}
{\LineFive{0}{2}{2}{1}{2}}
{\LineFive{0}{2}{0}{0}{1}}
{\LineFive{0}{2}{1}{1}{2}}
{\LineFive{0}{1}{2}{0}{0}}, \qquad b^{(2)} = \MatrixFive
{\LineFive{2}{1}{1}{1}{1}}
{\LineFive{0}{1}{2}{1}{1}}
{\LineFive{0}{1}{2}{2}{2}}
{\LineFive{0}{0}{1}{0}{0}}
{\LineFive{0}{1}{1}{1}{0}}.$$
Notice that $(b^{(1)})^3 = (1+\varphi_1^4)^3 = 1+ \varphi_1^{12} =
1+ \varphi_1^2 = 1+\varphi_2 = b^{(2)}$, and $(b^{(2)})^3 =
(1+\varphi_1^2)^3 = 1+\varphi_1^6 = b^{(1)}$.

Rather than ranging over  all $u \in \mul{\F_{243}}/\mul{\F_3}$,
we take a generator.  One can check that $t^{121} = 1$ in
$\F_{243}$, while $t^{11} = t^3-t^2+t \neq 1$. In our basis,
multiplication by $t$ corresponds to the matrix
$$\theta = \MatrixFive
{\LineFive{0}{0}{0}{0}{1}}
{\LineFive{1}{0}{0}{0}{1}}
{\LineFive{0}{1}{0}{0}{0}}
{\LineFive{0}{0}{1}{0}{0}}
{\LineFive{0}{0}{0}{1}{0}}.$$

For $i = 1,2$ and $j = 0,\dots,120$, let $b^{(i)}_{j} = \theta^j
b^{(i)}\theta^{-j}$, and set $\Omega^{(\s_i)} =
\set{b^{(i)}_0,\dots,b^{(i)}_{120}}$. Let $\bar{\Omega}^{(i)}$
denote the union of $\Omega^{(\s_i)}$ and its inverses. Let
$\hat{\Omega}^{(i)}$ denote the products $b^{(i)}_{j_1}\ldots
b^{(i)}_{j_k}$, $k = 1,\dots,4$, which can be completed to a
product of $5$ matrices from $\Omega^{(i)}$ which equals $1$.
These products correspond to the
$\binomq{5}{1}{3}+\binomq{5}{2}{3}+\binomq{5}{3}{3}+\binomq{5}{4}{3}
= 121+1210+1210+121 = 2662$ vector subspaces of $\F_3^5$, and
$\hat{\Omega}^{(i)}$ constitute symmetric sets of generators for
$\PSL[5](\F_3)$.

\begin{Cor}
Let $b^{(1)}$ and $\theta$ be the $5\times 5$ matrices over $\F_3$
given above, and let $\bar{\Omega}^{(1)}$, $\hat{\Omega}^{(1)}$ be
the sets of $242$ and $2662$ matrices, respectively, defined
above. Let $\bar{\Omega}^{(2)}$ and $\hat{\Omega}^{(2)}$ be the
sets of matrices obtained by raising every matrix in
$\bar{\Omega}^{(1)}$, respectively $\hat{\Omega}^{(1)}$, to the
third power.

Then $\bar{\Omega}^{(1)}$ and $\bar{\Omega}^{(2)}$ are two
symmetric sets of $242$ generators of $\PSL[5](\F_3)$, which
define isospectral non-isomorphic Cayley complexes. Likewise
$\hat{\Omega}^{(1)}$ and $\hat{\Omega}^{(2)}$ are symmetric sets
of $2662$ generators with the same properties.
\end{Cor}

Notice that $\PSL[5](\F_3)$ is a group of size $\approx 2.3\cdot
10^{11}$.  Therefore, storing its Cayley graphs in a computer is a
difficult task, and computing the eigenvalues directly is nearly
impossible.

More generally, we can carry out a similar construction whenever
$q$ and $d$ are odd and co-prime, and $e = 1$. Following the above
mentioned algorithm, take $\beta = 1$, then $y(\lam) = (1+\lam)^d
- 1$. Take $\alpha = -2$, then $\gamma = y(\alpha) = -2$ and the
minimal polynomials are $p(\lam) = g(\lam) = \lam+2$. With $s = 1$
we then have $\bar{L} = \F_q$ and $x = -2$. Now let $\varphi_1$
denote the matrix representing the Frobenius automorphism (in any
chosen basis of $\F_{q^d}$ over $\F_q$). For every $i =
1,\dots,d-1$, let $z_i = (1+\beta x) \varphi_1^{i} =  -
\varphi_1^{i}$, and $b^{(i)} = 1 - z_i^{-1} = 1 + \varphi_1^{-i}$.
Finally fix a generator $\theta$ of $\mul{\F_{q^d}}/\mul{\F_q}$,
and define the sets $\bar{\Omega}^{(i)}$ and $\hat{\Omega}^{(i)}$
as above. Notice that $(b^{(i)})^q = 1 + \varphi_1^{- q i} = b^{(q
i)}$ where the upper index here is modulo $d$, and so the elements
of $\Omega^{(qi)}$ are obtained by exponentiating the elements of
$\Omega{(i)}$ to the power $q$.

\begin{Cor}
Let $q,d$ be odd and co-prime, and let $m$ denote the order of $q$
in $\mul{\Z/d\Z}/\sg{\pm 1}$. Let $\varphi_1$ denote the Frobenius
automorphism of $\F_{q^d}/\F_q$, and let $\theta$ be an element of
$\F_{q^d}$ which generates $\mul{\F_{q^d}}/\mul{\F_q}$. Let
$\bar{\Omega}^{(1)}$ be the set of $2\frac{q^d-1}{q-1}$ elements
of $\PSL[d](\F_q)$ defined as above. Let $\bar{\Omega}^{(q^i)}$
denote the set obtained from $\hat{\Omega}^{(1)}$ by element-wise
exponentiation to the power $q^i$, for $i = 1,\dots,m-1$.

Then the Cayley graphs of $\PSL[d](\F_q)$ with respect to
$\bar{\Omega}^{(q^i)}$ are isospectral (even as complexes) and
non-isomorphic.
\end{Cor}
A similar result holds for the sets $\hat{\Omega}^{(q^i)}$.

\def\US{A.~Lubotzky, B.~Samuels and U.~Vishne}
\def\BIBPapI{\US, {\it Ramanujan complexes of type $\tilde{A_d}$},
Israel J. of Math., to appear.}
\def\BIBPapII{\US, {\it Explicit constructions of Ramanujan complexes},
European J. of Combinatorics, to appear.}
\def\BIBPapIII{\US, {\it Division algebras and non-commensurable isospectral manifolds},
preprint.}
\def\BIBPapIV{\US, {\it Isospectral Cayley graphs of some simple groups},
preprint.}

\end{document}